\documentclass[12pt]{article}

\usepackage{latexsym}
\usepackage{graphicx}
\usepackage{framed}
\usepackage{threeparttable}
\usepackage{amssymb}
\usepackage{amsmath}
\usepackage{multirow}
\usepackage{hyperref}

\begin{document}

\newcommand{\noi}{\noindent}
\newcommand{\nn}{\nonumber}
\newcommand{\bd}{\begin{displaymath}}
\newcommand{\ed}{\end{displaymath}}
\newcommand{\bp}{\underline{\bf Proof}:\ }
\newcommand{\ep}{{\hfill $\Box$}\\ }
\newtheorem{1}{LEMMA}[section]
\newtheorem{2}{THEOREM}[section]
\newtheorem{3}{COROLLARY}[section]
\newtheorem{4}{PROPOSITION}[section]
\newtheorem{5}{REMARK}[section]
\newtheorem{20}{OBSERVATION}[section]
\newtheorem{10}{DEFINITION}[section]
\newtheorem{30}{RESULTS}[section]
\newtheorem{40}{CLAIM}[section]
\newtheorem{50}{ASSUMPTION}[section]
\newtheorem{60}{EXAMPLE}[section]
\newtheorem{70}{ALGORITHM}[section]
\newtheorem{80}{PROBLEM}
\newcommand{\be}{\begin{equation}}
\newcommand{\ee}{\end{equation}}
\newcommand{\ba}{\begin{array}}
\newcommand{\ea}{\end{array}}
\newcommand{\bea}{\begin{eqnarray}}
\newcommand{\eea}{\end{eqnarray}}
\newcommand{\bqn}{\begin{eqnarray*}}
\newcommand{\eqn}{\end{eqnarray*}}


\newcommand{\e} { \ = \ }
\newcommand{\leqs}{ \ \leq \ }
\newcommand{\geqs}{ \ \geq \ }
\def\theequation{\thesection.\arabic{equation}}
\def\bReff#1{{\bRm
(\bRef{#1})}}
\newcommand{\eps}{\varepsilon}
\newcommand{\sgn}{\operatorname{sgn}}
\newcommand{\sign}{\operatorname{sign}}
\newcommand{\Vol}{\operatorname{Vol}}
\newcommand{\Var}{\operatorname{Var}}
\newcommand{\Cov}{\operatorname{Cov}}
\newcommand{\vol}{\operatorname{vol}}
\newcommand{\var}{\operatorname{var}}
\newcommand{\cov}{\operatorname{cov}}
\renewcommand{\Re}{\operatorname{Re}}
\renewcommand{\Im}{\operatorname{Im}}
\newcommand{\bE}{{\mathbb E}}
\newcommand{\bR}{\mathbb{R}}
\newcommand{\bN}{{\mathbb N}}
\newcommand{\bC}{\mathbb{C}}
\newcommand{\bF}{\mathbb{F}}
\newcommand{\bQ}{{\mathbb Q}}
\newcommand{\bZ}{{\mathbb Z}}
\newcommand{\cA}{{\mathcal A}}
\newcommand{\cB}{{\mathcal B}}
\newcommand{\cC}{{\mathcal C}}
\newcommand{\cD}{{\mathcal D}}
\newcommand{\cE}{{\mathcal E}}
\newcommand{\cF}{{\mathcal F}}
\newcommand{\cG}{{\mathcal G}}
\newcommand{\cH}{{\mathcal H}}
\newcommand{\cI}{{\mathcal I}}
\newcommand{\cJ}{{\mathcal J}}
\newcommand{\cK}{{\mathcal K}}
\newcommand{\cL}{{\mathcal L}}
\newcommand{\cM}{{\mathcal M}}
\newcommand{\cN}{{\mathcal N}}
\newcommand{\cO}{{\mathcal O}}
\newcommand{\cP}{{\mathcal P}}
\newcommand{\cQ}{{\mathcal Q}}
\newcommand{\cR}{{\mathcal R}}
\newcommand{\cS}{{\mathcal S}}
\newcommand{\cT}{{\mathcal T}}
\newcommand{\cU}{{\mathcal U}}
\newcommand{\cV}{{\mathcal V}}
\newcommand{\cW}{{\mathcal W}}
\newcommand{\cX}{{\mathcal X}}
\newcommand{\cY}{{\mathcal Y}}
\newcommand{\cZ}{{\mathcal Z}}
\newcommand{\bx}{{\mathbf x}}
\newcommand{\by}{{\mathbf y}}
\newcommand{\bz}{{\mathbf z}}
\newcommand{\bu}{{\mathbf u}}
\newcommand{\bv}{{\mathbf v}}
\newcommand{\bw}{{\mathbf w}}

\newcommand{\bbf}{{\mathbf f}}
\newcommand{\bbH}{{\mathbf H}}
\newcommand{\bzero}{{\mathbf 0}}
\newcommand{\bq}{{\mathbf q}}
\newcommand{\bba}{{\mathbf a}}
\newcommand{\bbb}{{\mathbf b}}
\newcommand{\bbc}{{\mathbf c}}


\title{A continuation method for tensor complementarity problems}
\author{ 
Lixing Han \\
Department of Mathematics \\
University of Michigan-Flint \\
Flint, MI 48502, USA \\
Email: \texttt{lxhan@umflint.edu}
}
 \date{March 1, 2018}
\maketitle

 \baselineskip 0.2 in

\begin{abstract}  
 We introduce a Kojima-Megiddo-Mizuno type continuation method for solving tensor complementarity problems. We show that there exists a bounded continuation trajectory when the tensor is strictly semi-positive and any limit point tracing the trajectory gives a solution of the tensor complementarity problem. Moreover, when the tensor is strong strictly semi-positive, tracing the trajectory will converge to the unique solution. Some numerical results are given to illustrate the effectiveness of the method.    
\end{abstract}



\newpage

\section{Introduction}
\label{Intro}
\setcounter{equation}{0}
Let $\bR^n$ denote the $n$-dimensional Euclidean space, $\bR_{+}^n = \{ \bx \in \bR^n : \bx \geq \bzero \}$ and $\bR_{++}^n = \{ \bx \in \bR^n : \bx > \bzero \}$. Let $\bbf: \bR^n \rightarrow \bR^n$ be a C$^1$ mapping, i.e., $\bbf$ is continuously differentiable. We consider the complementarity problem (CP) with the mapping $\bbf$: Find $\bx \in \bR^n$ such that
 \be
\label{cp}
\bx \geq \bzero, \ \ \ \bbf(\bx) \geq \bzero, \ \ \ \langle \bx, \bbf(\bx) \rangle = 0,
\ee
where $\langle \cdot,\cdot \rangle$ denotes the inner product on $\bR^n$. 
 
Let $\bR^{[m,n]}$ denote the set of all $m$th-order, $n$-dimensional real tensors. When $\bbf(\bx)=\cA \bx^{m-1} + \bq$, where $\cA \in \bR^{[m,n]}$,   $\bq \in \bR^n$, and $\cA \bx^{m-1}$ denotes the column vector  whose $i$th entry is
$$
(\cA \bx^{m-1})_i = \sum_{i_2, \cdots, i_m =1}^{n} A_{i  i_{2}  \cdots  i_m}x_{i_2} \cdots x_{i_m},   \ \ \ i=1,2, \ldots, n,
$$ 
the CP (\ref{cp}) becomes the so-called tensor complementarity problem (TCP), denoted by TCP($\cA, \bq$): Find $\bx \in \bR^n$ such that
\be
\label{tcp}
\bx \geq \bzero, \ \ \  \cA \bx^{m-1} + \bq \geq \bzero, \ \  \      \langle \bx,  \cA \bx^{m-1} + \bq \rangle = 0.
\ee

When $m=2$, the tensor complementarity problems reduce to the well studied linear complementarity problems \cite{CPS09}. When $m \geq 3$,  they form a nontrivial  class of nonlinear complementarity problems, which have received considerable attention recently (\cite{BHW16, CQS18, CQW16, DLQ15, GZH18, HQ17, LLV17, LQX15, SQ15, SQ16, SQ17, SY16, WCW18, WHB16}).  They have found applications in several areas, including nonlinear compressed sensing and game theory (\cite{LQX15, HQ17}). 

Structured tensors play an important role in studying theoretical properties of tensor complementarity problems.  Various structured tensors have been proposed and investigated in the literature, such as $P$ tensors and $P_0$ tensors, strong $P$ and $P_0$ tensors, $B$ tensors, $Q$ tensors, $R$ tensors, $S$ tensors, semi-positive tensors and strictly semi-positive tensors,  strong  semi-positive tensors and strong strictly semi-positive tensors, etc.. Of particular interest in this paper are the strictly semi-positive tensors and strong strictly semi-positive tensors as defined in the following definition.

\begin{10}  [\cite{SQ15, LLV17}]
\label{def1}
{\rm
Let $\cA \in \bR^{[m,n]}$. Then $\cA$ is called \\
(a) a strictly semi-positive tensor if for every $\bx \ne \bzero$ in $ \bR_{+}^{n}$, there is an index $i: 1 \leq 1 \leq n$ such that $x_i >0$ and
$(\cA \bx^{m-1} )_i > 0  $;   \\
(b) a strong strictly semi-positive if $\bbf(\bx) = \cA \bx^{m-1} +\bq$ is a $P$ function in $\bR_{+}^n$ for any $\bq \in \bR^n$, i.e., if for any distinct $\bx \in \bR_{+}^n$ and
  $\by \in \bR_{+}^n$, 
$$
\max_{1\leq i \leq n} (x_i-y_i) (f_i(\bx) - f_i(\by)) >  0. 
$$
}
\end{10}
  
It is easy to see that strong strictly semi-positiveness implies strictly semi-positiveness but the converse is not true (\cite{LLV17}). Song, Qi, and Yu \cite{SQ16, SQ17-0, SY16} proved the existence and boundedness of solutions of TCPs for strictly semi-positive tensors. Liu, Li, and Vong \cite{LLV17} proved that a TCP possesses the global uniqueness and solvability property if the tensor is strong strictly semi-positive, extending a similar result of Bai, Huang, and Wang \cite{BHW16} regarding strong $P$ tensors. We summarize their results in the following theorem.  

\begin{2} [\cite{SQ17-0, SY16, LLV17}]
\label{existence}
Let $\cA \in\bR^{[m,n]}$. Then for any $\bq \in \bR^n$, \\
(a) if $\cA$ is strictly semi-positive, the TCP($\cA, \bq$) has a nonempty compact solution set; \\
(b) if $\cA$ is strong strictly semi-positive, the TCP($\cA, \bq$) has a unique solution. 
\end{2}

Some algorithms for solving tensor complementarity problems have been proposed recently. Luo, Qi, Xiu \cite{LQX15} proposed a method for finding the sparsest solution to a TCP with a $Z$-tensor by reformulating the TCP as an equivalent polynomial programming problem.  Xie, Li, and Xu \cite{XLX17} proposed an iterative method for finding the least solution to a TCP.  Liu, Li, and Vong \cite{LLV17} proposed a modulus-based nonsmooth Newton's method for solving TCPs. Huang and Qi \cite{HQ17} proposed a smoothing type algorithm. 

Continuation methods form an important class of methods for solving linear or nonlinear complementarity problems (see \cite{KMM93, KMN91, KMN89, XD08, ZL01}). Recently homotopy continuation methods have been successfully developed for solving tensor eigenvalue problems and multilinear systems (see \cite{CHZ16, Han17}). Motived by this,  we consider using a homotopy continuation method for solving TCPs in this paper.  The convergence analyses of the continuation methods for general CPs in the literature impose certain conditions on the mapping $\bbf$ to ensure the existence and boundedness of a continuation trajectory (\cite{KMM93, KMN91, KMN89, XD08, ZL01}). Unfortunately, those analyses do not cover the TCP case when $\cA$ a  strictly semi-positive tensor, a condition that guarantees the existence and boundedness of solutions of the TCP.  It is therefore desirable that a continuation method for TCPs produces a bounded continuation trajectory when $\cA$ is strictly semi-positive. We will prove that this is indeed true in this paper.  We will also obtain some stronger results when $\cA$ is strong strictly semi-positive. We will implement the continuation method for TCPs with a strong strictly semi-positive tensor using an Euler-Newton predictor-corrector approach and provide some numerical results.  

This paper is organized as follows. In Section 2, we  introduce a Kojima-Megiddo-Mizuno type continuation method for solving TCPs and prove the existence of a bounded continuation path when the tensor $\cA$ is strictly semi-positive. Stronger results are proved when $\cA$ is strong strictly semi-positive. In Section 3, we present an implementation of the continuation method when $\cA$ is strong strictly semi-positive and some numerical results.  Some final remarks are given in Section 4.  

\section{A continuation method}
\label{method}
\setcounter{equation}{0}

  To solve a CP,  typically a continuation method first reformulates it as an equivalent problem. Here we use one of the most frequently used equivalent reformulations in the literature for the TCP (\ref{tcp}):  Find a solution $(\bx,\by) \in \bR^{2n}$ such that 
\be
\label{ecp}
  \left [ \begin{array}{c} 
                                      X \by \\
                                    \by-(\cA \bx^{m-1} + \bq)
   \end{array}
  \right ]   = \bzero,  \ \ \ (\bx,\by) \geq \bzero, 
\ee 
where $X={\rm diag}(\bx) $ is the diagonal matrix formed by the components of $\bx$. We then choose vectors $\bba \geq \bzero$ and $\bbb > \bzero$ from $\bR^n$ and define a Kojima-Megiddo-Mizuno (\cite{KMM93}) type homotopy mapping $\bbH: \bR^{2n} \times [0,1]
\rightarrow \bR^{2n} $ by 
\be
\label{homoH}
\bbH(\bx,\by,t) = \left [ \begin{array}{c} 
                                      X\by  - t \bba\\
                                    \by - (1-t)  (\cA \bx^{m-1} + \bq ) - t \bbb
   \end{array}
  \right ]. 
\ee

Starting with $t=1$ and $(\bx^0,\by^0)=(B^{-1}\bba, \bbb)$, where $B = {\rm diag}(\bbb)$, the continuation method follows a path from $t=1$ to $t=0$ by solving the system 
\be
\label{homosystem}
\bbH(\bx,\by,t)=0, \ \ \ (\bx,\by) \geq \bzero. 
\ee

Let $\bz = (\bx, \by)$. Denote $A={\rm diag}(\bba), B={\rm diag}(\bbb)$, $Y={\rm diag} (\by)$. 
The partial derivatives matrices $D_{\bz} \bbH(\bx,\by, t)$ and $D_{t} \bbH(\bx,\by, t)$ of the homotopy $\bbH(\bx, \by, t)$ play an important role in solving the system (\ref{homosystem}). To compute $D_{\bz} \bbH(\bx,\by, t)$, we need  the so-called semi-symmetric tensor $\hat{\cA}=(\hat{A}_{i_1, i_2, \ldots, i_m})$   (\cite{NQ15}) defined by
\be
\label{psymm}
\hat{A}_{i_1 i_2 \ldots i_m} = \frac{1}{p} \sum_{k=1}^{p} A_{i_1 i_{2}^{(k)}\ldots   i_{m}^{(k)}},
\ee
 where the sum is over all the $p$ different permutations $i_{2}^{(k)}, \ldots,  i_{m}^{(k)}$ of $i_2, \ldots, i_m$.   The partial derivatives matrix of $\cA \bx^{m-1}$ with respect to $\bx$ is 
\be
\label{df}
D_\bx \cA \bx^{m-1} = (m-1) \hat{\cA} \bx^{m-2}. 
\ee
 Therefore, the partial derivatives of $\bbH$ with respect to  $\bz$ and $t$ are:
\be
\label{jacobian}
D_{\bz} \bbH(\bx, \by,t) = \left [  \begin{array}{cc} 
Y  & X \\
-(1-t)(m-1) \hat{\cA} \bx^{m-2} & I_n
\end{array}
 \right ]
\ee
and 
\be
\label{partialT}
D_{t} \bbH(\bx, \by,t) = \left [  \begin{array}{c} 
-\bba \\
 \cA \bx^{m-1} + \bq - \bbb
\end{array}
 \right ]
\ee
respectively, where $I_n \in \bR^{n \times n}$ is the identity matrix.

When $t=1$, the matrix 
\[
D_{\bz} \bbH(\bx^0, \by^0,1) = \left [  \begin{array}{cc} 
B  & B^{-1} A \\
0 & I
\end{array}
 \right ]
\]
is nonsingular. By the Implicit Function Theorem, there is $\delta_1 \in (0,1)$ such that the system (\ref{homosystem}) has a unique solution $(\bx(t), \by(t))$ for each $t \in (1-\delta_1,1]$ such that $\bx(1) = \bx^0, \by(1)=\by^0$, and $\bx(t)$ and $\by(t)$ are smooth functions of $t$. Thus, the homotopy (\ref{homosystem}) has a unique smooth trajectory emanated from $(\bx^0, \by^0,0)$ for $t \in [1-\delta,1]$, where we choose $\delta \in (0, \delta_1)$. Denote this trajectory by 
\be
\label{initialT}
T_{\delta} = \{ (\bx (t), \by(t), t):  t \in [1-\delta,1] \}.
\ee

Under rather mild conditions, we can show the existence of a continuation path that contains $T_{\delta}$. We summarize the result in the following theorem, whose proof can be found in \cite{KMM93} for a general CP.  

\begin{2} [\cite{KMM93}]
\label{thm1} 
Let $\bba \in \bR_{+}^n$ be fixed. Then, 
for almost every $\bbb \in \bR_{++}^n$, starting from $(\bx^0, \by^0,1)$, the homotopy system (\ref{homosystem}) yields a trajectory 
\be
\label{trajectory}
T = \{ (\bu(s), \bv(s), t(s))  
\in \bR_{+}^{n} \times \bR_{+}^n \times (0,1]:  0< s \leq1 \},
\ee
which contains $T_{\delta}$. Here $\bu: (0,1] \to \bR_{+}^n$, $\bv: (0,1] \to \bR_{+}^n$, and $t: (0,1] \to (0,1]$ are piecewise $C^1$ mappings, and $T$ is a 1-dimensional manifold that is homeomorphic to $(0,1]$.  

If $T$ is bounded, then $\lim_{s \to 0} t(s) = 0$. 

If $\bba \in \bR_{++}^n$, then $\bu$, $\bv$, and $t$ are $C^1$ mappings.

\end{2}

\begin{5}
\label{remark1}
{\rm
When $\bba \in \bR_{++}^n$, the trajectory  $T$ in Theorem \ref{thm1} is a 1-dimensional smooth manifold. It is known that a 1-dimensional smooth manifold is diffeomorphic to a unit circle or a unit interval (see, for example, \cite{Naber80}).   Since the matrix $D_{\bz} \bbH(\bx^0, \by^0,1)$ is nonsingular, $T$ is not diffeomorphic to a unit circle. Hence, $T$ is  diffeomorphic to $(0,1]$. 
}
\end{5}

To ensure the boundedness of a trajectory in a continuation method  for a general CP (\ref{cp}), certain conditions need to be imposed on the mapping $\bbf$. Among several such conditions proposed in the literature, the following two conditions have been frequently used in theoretical studies and practice. 

\ \\
{\bf Condition 1.} (\cite{KMM93, ZL01})  \\
(a) $\bbf$ is monotone on $\bR_{+}^n$, i.e., $\langle \bx - \by,  \bbf(\bx) - \bbf(\by)  \rangle \geq 0$ for any $\bx$ and $\by$ in $\bR_{+}^n$.   \\
(b) There exists a strictly feasible point $(\bar{\bx}, \bar{\by})$ such that 
$$\bar{\bx} >\bzero \ {\rm and} \  \bar{\by} = \bbf(\bar{\bx}) >\bzero.$$

\ \\
{\bf Condition 2.} (\cite{KMN91, KMN89, ZL01}) \\
(a) $\bbf$ is a $P_0$ function in $\bR_{+}^n$, i.e., for any distinct $\bx$ and $\by$ in $\bR_{+}^n$, 
$$
\max_{x_i \ne y_i}  (x_i-y_i) (f_i(\bx) - f_i(\by) ) \geq 0. 
$$
 (b)  There exists a strictly feasible point $(\bar{\bx}, \bar{\by})$. \\
(c) The set 
$$
\bbf^{-1}(D) = \{\bz=(\bx, \by) \in \bR_{+}^{2n}: \bbf(\bz) \in D   \}  
$$
is bounded for every compact subset of $D$ of $\bR_{+}^{n} \times B_{++}(\bbf)$, where 
$$
B_{++} (\bbf) = \{ \bu= \by - \bbf (\bx) \ {\rm for \ all} \ (\bx, \by) > \bzero   \}.
$$

The condition imposed in \cite{XD08} does not require that $\bbf$ is a $P_0$ mapping or a monotone mapping. However, they require the following:

\ \\
{\bf Condition 3.} (\cite{XD08}) \\
(a) $\bbf$ is three times continuously differentiable. \\
 (b)  For any $\{\bx^k\} \subset \bR_{+}^n $, as $k \to \infty$ and $\|\bx^k\| \to \infty$, $\bbf(\bx^k) >\bzero $ when $k > K_0$ for some $K_0>0$.   \\
(c) The set $S_{+} = \{ (\bx,\by) \geq \bzero: \by = \bbf(\bx) \}$ is non-empty. \\

Unfortunately, there exist  strictly semi-positive tensors such that none of these conditions hold. As an example, we take Example 3.27 in \cite{LLV17} with  the following $\cA \in \bR^{[3,2]}$:
$$a_{111}=1, a_{121}=2, a_{122}=1, a_{222}=1, a_{211}=-1, a_{221}=-1,$$ 
and
$$
a_{ijk}=0,  \ {\rm for \ other} \ i,j,k.
$$
It is shown in \cite{LLV17} that this tensor is strictly positive semi-definite. However,  its corresponding mapping  $\bbf(\bx)=\cA \bx^{m-1} + \bq$ is not a $P_0$ function or a monotone function for any $\bq \in \bR^n$.  
Moreover, as
\[
\cA \bx^2 = \left [ \begin{array}{c}
(x_1+x_2)^2 \\
x_{2}^2-x_1x_2-x_{1}^2
\end{array}
\right ],
\]
clearly $\bbf(\bx)=\cA \bx^{2}$ does not satisfy Condition (3.b).  

Thus,  the proofs of boundedness of the trajectory for a general CP in the literature do not apply to the TCP case when $\cA$ is strictly semi-positive. According to Theorem \ref{existence}, $TCP(\cA, \bq)$ has a nonempty compact solution set if $\cA$ is strictly semi-positive. Therefore, it is desirable to prove  that the trajectory (\ref{trajectory}) is bounded if the tensor $\cA$ is strictly semi-positive. As we will show in the following theorem, this is indeed true.  

\begin{2}
\label{main}
Under the conditions of Theorem \ref{thm1},  if $\cA \in \bR^{[m,n]}$ is a strictly semi-positive tensor,  then the trajectory  (\ref{trajectory})   
is bounded. As $s \to 0$, any limit point of $(\bu(s), \bv(s))$ tracing this trajectory is a solution of the problem (\ref{ecp}).  
\end{2}
\bp
Clearly, for $t(s) \in [1-\delta, 1]$, $T=T_{\delta}$ is bounded. Thus, we only need to prove that the trajectory
$$
T = \{ (\bu(s), \bv(s), t(s)):  0 < s \leq 1 \}
$$ 
is bounded when $ t(s) \in (0,1-\delta]$.

Since $\cA$ is strictly semi-positive, as in \cite{SY16}, we define the following quantity
$$
\beta(\cA)= \min_{\begin{array}{c} 
\bx \geq \bzero \\ \|\bx \|_{\infty}=1
\end{array}}  \max_{1\leq i \leq n} x_i (\cA \bx^{m-1})_{i}.
$$
Then $\beta(\cA)>0$.  

Note that 
$$
v(s)_i = (1-t(s)) (\cA \bu(s)^{m-1}+\bq )_i  + t(s) b_i,
$$
and
$$
u(s)_iv(s)_i=t(s) a_i.
$$
Thus,
$$
t(s) a_i = (1-t(s)) u(s)_i (\cA \bu(s)^{m-1})_i   + u(s)_i ((1-t(s))  q_i + t(s)  b_i).
$$
By the definition of $\beta(\cA)$, we have
\bqn 
(1-t(s)) \|\bu(s)\|_{\infty}^m \beta (\cA) & \leq & (1-t(s)) \max_{i} u(s)_i (\cA \bu(s)^{m-1})_i  \\
      & \leq & \max_i  t(s) a_i   + \max_i u(s)_i ( -(1-t(s)) q_i - t(s) b_i )   \\
     & \leq & \max_i a_i + \|\bu(s)\|_{\infty} \|- (1-t(s))\bq - t(s) \bbb \|_{\infty}  \\
   & \leq & \max_i a_i + \|\bu(s)\|_{\infty} ( \|\bq\|_{\infty} + \| \bbb \|_{\infty})
\eqn 
Therefore,
$$
\|\bu(s)\|_{\infty}^{m-1} \leq \max  \left \{1 , \frac{ \max_{i} a_i  + \|\bq\|_{\infty} + \| \bbb \|_{\infty}  }{ (1-t(s)) \beta(\cA) }  \right   \}.
$$
This implies that $\bu(s)$ is bounded when $t(s) \in (0,1-\delta]$. Thus, it is bounded for $ s \in (0,1]$. The boundedness of $\bv(s)$ follows from the boundedness of $\bu(s)$.

The boundedness of $T$ implies that 
$ \lim_{s \to 0} t(s)=0$ by Theorem \ref{thm1}. Moreover, $(\bu(s),\bv(s))$ has limit points when $s \to 0$. Each limit point is a solution of the problem (\ref{ecp}).  
\ep

We now study the behavior of the trajectory $T$ in (\ref{trajectory})  when the tensor $\cA$ is strong strictly semi-positive. By Remark \ref{remark1} and Theorem \ref{main}, $T$ is  a bounded smooth 1-dimensional manifold that is diffeomorphic to $(0,1]$. On this trajectory, we have 
$$
\bbH(\bu(s),\bv(s), t(s)) = \bzero,
$$ 
for $s \in (0,1]$. Let $\bz= (\bu, \bv)$. Since $\bu$, $\bv$, and $t$ are smooth in $s$, differentiating this system gives
\be
\label{diffeq}
D_{\bz} \bbH \cdot \frac{d \bz}{d s} + D_t \bbH \cdot \frac{d t}{d s} =\bzero. 
\ee
Following (\ref{jacobian}), the partial derivatives matrix $D_{\bz} \bbH$ is of the form
\be
\label{jacobian2}
 D_{\bz} \bbH( \bu(s),\bv(s), t(s) )  = \left [  \begin{array}{cc} 
V  & U \\
-M & I_n
\end{array}
 \right ],
\ee
where $V={\rm diag} (\bv(s))$, $U={\rm diag} (\bu(s))$, and $M= (1-t(s))(m-1) \hat{\cA} \bu(s)^{m-2} $. Since $\cA$ is is strong strictly semi-positive, $\cA \bx^{m-1}$ is a $P$ function in $\bR_{+}^n$. This implies that its Jacobian matrix $D_{\bx}  \cA \bx^{m-1}$ defined in (\ref{df}) is a $P_0$ matrix for $\bx \in \bR_{+}^n$ (see, \cite{More96}). Now as $\bu(s) \in \bR_{++}^n$, $\bv(s) \in \bR_{++}^n$, and $t(s) \in (0,1]$, the product matrix $VU^{-1}$ is a diagonal matrix with positive diagonal entries and the matrix $M$ is a $P_0$ matrix. It follows that the partial derivatives matrix $D_{\bz} \bbH$ in (\ref{diffeq}) along the trajectory is nonsingular, as  its determinant
$$
{\rm det} (D_{\bz} \bbH( \bu(s),\bv(s), t(s) ) ) = {\rm det} ( V+MU  ) =  {\rm det} ( VU^{-1}+M  ) {\rm det} (U) \ne 0.
$$
This means that $dt/ds \ne 0$ for any $s \in (0,1]$, because otherwise the nonsingularity of  $D_{\bz} \bbH$ would imply $d \bz/ds =0$ for some $s $, contradicting that $T$ is diffeomorphic to $(0,1]$.  
Therefore, the trajectory
$T$ can be parametrized by the variable $t$: 
\be
\label{traject2}
T= \{ (\bx(t), \by(t), t): 0 < t \leq 1,    \} \subset \bR_{+}^n \times \bR_{+}^n \times (0,1], 
\ee
and we can trace $T$ by solving the following initial value problem
\be
\label{ivp}
D_{\bz} \bbH(\bz, t)  \frac{d \bz}{d t} =- D_t \bbH (\bz, t), \ \  \ \  \bz(1) =(\bx^0, \by^0),  
\ee 
where $\bz = (\bx, \by)$. 

We summarize the convergence results of the continuation method when $\cA$ is strong strictly semi-positive in the following theorem. 

\begin{2}
\label{ssp-unique}
Let $\bba \in \bR_{++}^n$. Suppose that $\cA \in \bR^{[m,n]}$ is a strong strictly semi-positive tensor and $\bq \in \bR^n$.  Then for almost every $\bbb \in \bR_{++}^n$, solving the system (\ref{homosystem}) yields a a smooth and bounded trajectory (\ref{traject2}). Moreover,   tracing this trajectory by solving the initial value problem (\ref{ivp})  converges to the unique solution $(\bx_*, \by_*)$ of the problem (\ref{ecp}) as $t \to 0$, in which $\bx_*$ is the unique solution of the TCP (\ref{tcp}). 
\end{2}

\begin{5}
\label{remark2}
{\rm
The analysis given before Theorem \ref{ssp-unique} is valid when the tensor$\cA$ is strictly semi-positive and  the function $\bbf (\bx) = \cA \bx^{m-1} +\bq$ is a $P_0$ function in $\bR_{+}^n$, i.e., 
if for any distinct $\bx \in \bR_{+}^n$ and
  $\by \in \bR_{+}^n$, 
$$
\max_{1\leq i \leq n} (x_i-y_i) (f_i(\bx) - f_i(\by)) \geq   0. 
$$
In this case, the solution set of the TCP (\ref{tcp}) (as well as the solution set of (\ref{ecp})) is nonempty and compact, but the solutions are not necessarily unique. Therefore, any limit point of the trajectory (\ref{traject2}) as $t \to 0$ is a solution of (\ref{ecp}).    
}
\end{5}

\section{Numerical results} 
\label{NumResults}
\setcounter{equation}{0}

We have implemented the continuation method described in the previous section when the tensor $\cA$ is strong strictly semi-positive. An Euler-Newton  predictor-corrector method with adaptive step sizes (see, for example, \cite{AG90})  is used to solve the initial value problem (\ref{ivp}). We summarize our implementation in the following algorithm.

\begin{70}
\label{algorithmA}

\ \\
\ \\
{\bf  Step 0.} (Initialization)  Choose positive vectors $\bba, \bbb \in \bR_{++}^{n}$.  Choose initial step size $\Delta t_0>0$, tolerances $\epsilon_1>0$ and $\epsilon_2>0$. Let 
$t_0=1$, $\bx_0=B^{-1} \bba$, $\by_0=\bbb$. Let $\bz = (\bx, \by)$ and $ \bz_0 = (\bx_0, \by_0)$. Set $k=0$. 
\ \\
{\bf Step 1.} (Find new $t_{k+1}$): Set $ t_{k+1} = t_k - \Delta t_k$. If $t_k>0$ and $t_{k+1} \leq 0$ for some $k$, then set $N=k$ and reset $t_{N+1} = 0$ and $\Delta t_N=t_N$.
\ \\
{\bf Step 2.} (Find a predictor using Euler's method) Compute the tangent vector $ \mathbf{ g} $ to $\bbH(\bz,t)=0$ at $t_k$ by solving the linear system
\begin{equation*}
D_{\bz} \bbH(\bz_k,t_k) \mathbf{g} = -D_t  \bbH(\bz_k,t_k)
\end{equation*}
for $\mathbf{g}$. Then compute the approximation 
$\tilde{\bz}$ to $\bz_{k+1}$ by
$$
\tilde{\bz} = \bz_k + \Delta t_k \mathbf{g}.
$$
\ \\
{\bf Step 3.} (Find a corrector using Newton's method)  Initialize $\bw_0 = \tilde{\bz}$. For $i = 0, 1, 2, \dots$, compute
$$
\bw_{i+1} = \bw_i - [D_{\bz} \bbH(\bw_i,t_{k+1})]^{\dagger} \bbH(\bw_i,t_{k+1})
$$
until $\|\bbH(\bw_{J},t_{k+1})\| \leq \epsilon_1$ if $k<N$ or $\|\bbH(\bw_{J},t_{k+1})\| \leq \epsilon_2$ if $k=N$, where $\dagger$ denotes the pseudo-inverse. Then let $\bz_{k+1} = \bw_{J}$. If $k=N$, we set $(\bx_*,\by_*)= \bz_{N+1} $ as the computed solution of problem (\ref{ecp}) and stop.
\ \\
{\bf Step 4.} (Adaptively update the step size $\Delta t_k$)  If more than three steps of Newton iterations were required to converge within the desired accuracy, then $\Delta t_{k+1}=0.5 \Delta t_k$. If $\Delta t_{k+1} \leq 10^{-6}$, set $\Delta t_{k+1}=10^{-6}$.  If two consecutive steps were not cut, then  $\Delta t_{k+1}=2 \Delta t_k$. If $\Delta t_{k+1} \geq 0.5$, set $\Delta t_{k+1}=0.5$. Otherwise, $\Delta t_{k+1}=\Delta t_k$. 
 Set $k=k+1$. Go to Step 1.  

\end{70}

We have coded Algorithm \ref{algorithmA} in MATLAB and done some numerical experiments on the following  examples. 

\begin{60} 
\label{example1} 
{\rm
  Let  $\cA \in \bR^{[3,2]}$ be defined by:
$$
a_{111}=1, a_{121}=1, a_{122}=-1, a_{222}=1,  a_{211}=-1, a_{221}=1, 
$$
and  $ a_{i_1i_2i_3}=0$  otherwise.  This tensor is given in \cite[Example 3.30]{LLV17}. It is strong strictly semi-positive. Different vectors $\bq \in \bR^2$ are used in our experiments.
}
\end{60}

\begin{60} 
\label{example2} 
{\rm
  Consider the tensor $\cA \in \bR^{[4,2]}$ defined by:
$$
a_{1111}=1, a_{1222}=-1, a_{1122}=1, a_{2222}=1,  a_{2111}=-1, a_{2211}=1, 
$$
and  $ a_{i_1i_2i_3i_4}=0$ otherwise.  This tensor is first given in \cite[Example 4.2]{BHW16}. It
 is a $P$ tensor, but not a strong $P$ tensor. However,   it is strong strictly semi-positive \cite{LLV17}. Different vectors $\bq \in \bR^2$ are used in our experiments.
}
\end{60}

\begin{60}
\label{example3}
{\rm
 Let   $\cA \in \bR^{[5,3]}$ be defined by $a_{kkkkk}=k$, for $k=1,2,3$, and $a_{i_1i_2i_3i_4i_5}=0$ otherwise. Clearly, this tensor is strong strictly semi-positive. Different vectors $\bq \in \bR^3$ are used in our experiments. 
}
\end{60}

\begin{60} 
\label{example4} 
{\rm
  Consider the tensor $\cA \in \bR^{[3,2]}$ defined by:
$$
a_{111}=1, a_{121}=2, a_{122}=1, a_{222}=1,  a_{211}=-1, a_{221}=-1, 
$$
and  $ a_{i_1i_2i_3}=0$ otherwise.  This tensor is given in \cite[Example 3.27]{LLV17}. It
 is strictly semi-positive,  but not strong strictly semi-positive. Nonetheless, the TCP (\ref{tcp}) with this tensor has a unique solution for any $\bq \in \bR^2$.    Different vectors $\bq \in \bR^2$ are used in our experiments.
}
\end{60}

Our experiments were done using MATLAB 2014b on a laptop computer with Intel Core i7-4600U at 2.10 GHz and 8 GB memory running Microsoft Windows 7.  The tensor toolbox of \cite{BK15} was used to compute tensor-vector products and to compute the semi-symmetric tensor $\hat{\cA}$.  We used $\bba=\bbb=[1,1,\ldots,1]^T$, $\Delta t_0=0.1$, $\epsilon_1=10^{-5}$, and $\epsilon_2=10^{-12}$ in Algorithm \ref{algorithmA}.

We now report the numerical results in Tables \ref{table1}, \ref{table2}, \ref{table3}, and \ref{table4}.  In these tables,  \texttt{itr} and \texttt{nwtitr} denote the number of prediction steps and the number of Newton iterations were used, respectively, \texttt{solution} denotes the solution of TCP (\ref{tcp}) found by Algorithm \ref{algorithmA}, and \texttt{residue} denotes the residue
\[
\left \|    \left [ \begin{array}{c} 
                                      X_* \by_* \\
                                    \by_*-(\cA \bx_{*}^{m-1} + \bq)
   \end{array}
  \right ]    \right  \|_2
\]
 at termination, where $X_* = {\rm diag} (\bx_*) $.

\begin{table}[htbp]
\begin{center}
\begin{tabular}{|c||c|c|c|c|}
\hline
$\bq$ & itr & nwtitr & solution & residue  \\ \hline
$[-5,-3]^T$  & 5 & 12 & $[2.1286,1.8792]^T$ 	&  $8.8805e-15$  \\ \hline
$[-5,3]^T$  & 5 & 14 & $[2.0582,0.4859]^T$ 	&  $2.6746e-23$  \\ \hline
$[5,3]^T$  & 5 & 12 & $[0,0]^T$ 	&  $4.2730e-13$  \\ \hline
$[0,3]^T$  & 5 & 36 & $[0,0]^T$ 	&  $6.8709e-13$  \\ \hline
$[2,-3]^T$  & 5 & 13 & $[0.3103,1.6113]^T$ 	&  $8.9179e-16$  \\ \hline
$[0,-5]^T$  & 5 & 12 & $[1.2430,2.0112]^T$ 	&  $2.1817e-15$  \\ \hline
\end{tabular}
\end{center}
\caption{Numerical Results for Example \ref{example1}}
\label{table1}
\end{table}

\begin{table}[htbp]
\begin{center}
\begin{tabular}{|c||c|c|c|c|}
\hline
$\bq$ & itr & nwtitr & solution & residue  \\ \hline
$[-5,-3]^T$  & 5 & 13 & $[1.6678,1.5096]^T$ 	&  $1.1586e-14$  \\ \hline
$[-5,3]^T$  & 5 & 14 & $[1.6714,0.5409]^T$ 	&  $1.9860e-15$  \\ \hline
$[5,3]^T$  & 5 & 12 & $[0,0]^T$ 	&  $1.0731e-18$  \\ \hline
$[0,3]^T$  & 5 & 33 & $[0,0]^T$ 	&  $5.2577e-13$  \\ \hline
$[2,-3]^T$  & 5 & 13 & $[0.3906,1.4167]^T$ 	&  $6.2804e-16$  \\ \hline
$[0,-5]^T$  & 5 & 13 & $[1.1143,1.6331]^T$ 	&  $3.8998e-15$  \\ \hline
\end{tabular}
\end{center}
\caption{Numerical Results for Example \ref{example2}}
\label{table2}
\end{table}

\begin{table}[htbp]
\begin{center}
\begin{tabular}{|c||c|c|c|c|}
\hline
$\bq$ & itr & nwtitr & solution & residue  \\ \hline
$[1,2,3]^T$  & 5 & 12 & $[0,0,0]^T$ 	&  $4.0969e-21$  \\ \hline
$[1,-2,3]^T$  & 5 & 12 & $[0,1,0]^T$ 	&  $7.6027e-23$  \\ \hline
$[-3,-2,-3]^T$  & 5 & 11 & $[1.3161,1,1]^T$ 	&  $3.6186e-15$  \\ \hline
$[3,3,3]^T$  & 5 & 12 & $[0,0,0]^T$ 	&  $4.9693e-23$  \\ \hline
$[-3,-1,-2]^T$  & 5 & 11 & $[1.3161,0.8409,0.9036]^T$ 	&  $3.6748e-15$  \\ \hline
$[0,-1,-2]^T$  & 5 & 31 & $[0,0.8409,0.9036]^T$ 	&  $7.1789e-13$  \\ \hline
\end{tabular}
\end{center}
\caption{Numerical Results for Example \ref{example3}}
\label{table3}
\end{table}

\begin{table}[htbp]
\begin{center}
\begin{tabular}{|c||c|c|c|c|}
\hline
$\bq$ & itr & nwtitr & solution & residue  \\ \hline
$[-5,-3]^T$  & 5 & 14 & $[0.3127,1.9233]^T$ 	&  $1.2942e-15$  \\ \hline
$[-5,3]^T$  & 5 & 11 & $[1.5513,0.6847]^T$ 	&  $1.7402e-14$  \\ \hline
$[5,3]^T$  & 5 & 13 & $[0,0]^T$ 	&  $6.0454e-26$  \\ \hline
$[0,3]^T$  & 5 & 36 & $[0,0]^T$ 	&  $6.3603e-13$  \\ \hline
$[2,-3]^T$  & 5 & 13 & $[0,1.7321]^T$ 	&  $9.9301e-16$  \\ \hline
$[0,-5]^T$  & 5 & 15 & $[0,2.2361]^T$ 	&  $1.9860e-15$  \\ \hline
\end{tabular}
\end{center}
\caption{Numerical Results for Example \ref{example4}}
\label{table4}
\end{table}

From these tables, we observe that Algorithm \ref{algorithmA} effectively computes the unique solution for each TCP in Examples \ref{example1}--\ref{example4}. The algorithm is also efficient in terms of the number of prediction steps \texttt{itr} and the number of Newton iterations \texttt{nwtitr}. We remark that the relatively large \texttt{nwtitr} in the cases when $\bq=[0,3]^T$ in Examples \ref{example1}, \ref{example2}, and \ref{example4}, and $\bq=[0, -1,-2]^T$ in Example \ref{example3}  is because more Newton iterations were used in the last step due to the singularity of the Jacobian matrix $D_{\bz} H(\bx_*, \by_*, 0)$ at the solution $(\bx_*, \by_*)$ of problem (\ref{ecp}).  Using a deflation method such as the one given in \cite{LVZ06} can improve the performance of Newton's method in such cases.

\section{Concluding Remarks}
\label{sec4}

We have introduced a continuation method for solving TCPs.  Under the assumption that the tensor is strictly semi-positive, we have proved the existence of a bounded continuation trajectory. This result is not covered by the theoretical results proved in the literature for general nonlinear complementarity problems. We have also proved that when the tensor is strong strictly semi-positive, tracing the trajectory will converge to the unique solution of the TCP. We have implemented the method for TCPs with strong strictly semi-positive tensors. Numerical results show the continuation method is promising for solving TCPs.

Various structured tensors have been introduced recently and they play an important role in studying theoretical properties of TCPs.   An interesting direction for  future research  is to investigate how to use a continuation method to solve TCPs with other types of structured tensors.


\end{document}